\documentclass[a4paper,10pt]{amsart}

\usepackage{upref}

\usepackage[colorlinks,backref]{hyperref}                   

\usepackage{amssymb,amsmath}

\newtheorem{theorem}{Theorem}[section]

\newtheorem{definition}{Definition}[section]

\theoremstyle{remark} 
\theoremstyle{definition} 

\numberwithin{equation}{section}

\newif\ifcomment \commentfalse \def\commentON{\commenttrue}
 \long\outer\def\BC#1\EC{\ifcomment \sloppy \par
\# \ldots\dotfill {\em #1} \dotfill \# \par \fi } \commentON

\newcommand{\remove}[1]{}

\newcommand{\eps}{\ensuremath{\varepsilon}}

\newcommand{\R}{{{\mathbb{R}}}}

\newcommand{\Div}{\mathrm{div}}
\newcommand{\dxdt}{\, dxdt}
\newcommand{\dyds}{\, dyds}
\newcommand{\dxi}{\, d\xi}

\newcommand{\dx}{\, dx}

\newcommand{\ds}{\, ds}
\newcommand{\dt}{\, dt}
\newcommand{\dy}{\, dy}

\newcommand{\dz}{\, dz}
\newcommand{\dtau}{\, d\tau}

\newcommand{\fourint}{\int\!\!\!\!\int\!\!\!\!\int\!\!\!\!\int}
\newcommand{\threeint}{\int\!\!\!\!\int\!\!\!\!\int}

\newcommand{\ps}{\ensuremath{\partial_s}}
\newcommand{\pt}{\partial_t}
\newcommand{\pxi}{\partial_{x_i}}

\newcommand{\pyi}{\partial_{y_i}}
\newcommand{\pyj}{\partial_{y_j}}
\newcommand{\pyiyj}{\ensuremath{\partial_{y_iy_j}^2}}
\newcommand{\pxixj}{\ensuremath{\partial_{x_ix_j}^2}}
\newcommand{\pxiyj}{\ensuremath{\partial_{x_iy_j}^2}}

\newenvironment{Definitions}
{%

\begin{enumerate}}%
{\end{enumerate}
}

\newcommand{\sgn}{{\rm sgn}\, }

\newcommand{\Levy}{\ensuremath{\mathcal{L}}}

\newcommand{\abs}[1]{\left| #1 \right|}
\newcommand{\norm}[1]{\ensuremath{\left\|#1\right\|}}
\newcommand{\Set}[1]{\ensuremath{\left\{#1\right\}}}

\newcommand{\weak}{\rightharpoonup}

\def\dsp{\displaystyle}

\begin{document}

\title[Fractional degenerate parabolic equations]
{Stability of entropy solutions for L\'evy
mixed hyperbolic-parabolic equations}

\author[K. H. Karlsen]{Kenneth H. Karlsen}
\address[Kenneth H. Karlsen]{\newline
         Centre of Mathematics for Applications \newline
         University of Oslo\newline
         P.O. Box 1053, Blindern\newline
         N--0316 Oslo, Norway}
\email[]{kennethk@math.uio.no}
\urladdr{http://folk.uio.no/kennethk}

\author[{S. Ulusoy}]{Suleyman Ulusoy}
\address[ S\"{u}leyman Ulusoy]{\newline
	Centre of Mathematics for Applications (CMA)\newline
	Department of Mathematics\newline
	University of Oslo\newline
	P.O. Box 1053, Blindern\\ N--0316 Oslo, Norway}
\email{suleymau@cma.uio.no}
\urladdr{http://folk.uio.no/suleymau/}




\thanks{The work was supported by the Research Council of Norway through
an Outstanding Young Investigators Award of K. H. Karlsen. 
This article was written as
part of the the international research program on 
Nonlinear Partial Differential Equations at the Centre for
Advanced Study at the Norwegian Academy of Science
and Letters in Oslo during the academic year 2008--09.}

\date{\today}

\begin{abstract}
We analyze entropy solutions for a class of L\'evy
mixed hyperbolic-parabolic equations
containing a non-local (or fractional) diffusion operator
originating from a pure jump L\'evy process.
For these solutions we establish uniqueness ($L^1$ contraction property)
and continuous dependence results.
\end{abstract}

\maketitle

\tableofcontents

\section{Introduction}
\label{sec:intro}

The subject of this paper is uniqueness and stability
results for properly defined entropy solutions of
mixed hyperbolic-parabolic quasilinear equations appended with
a nonlocal (fractional) diffusion operator. 
These equations take the form
\begin{equation}\label{eq1}
	\partial_t u + \Div f(u) = \Div (a(u) \nabla u)
	+ \Levy[u],
\end{equation}
where $u=u(t,x)$ is the unknown, $(t,x) \in 
Q_T := (0,T) \times \R^d$, $d\ge 1$, and $T>0$ is a
fixed final time. The operator $\Levy$ is the generator 
of a symmetric positivity preserving
pure jump L\'evy semigroup $e^{t\Levy}$ on $L^1(\R^d)$.

Equation \eqref{eq1} is subject to initial data
\begin{equation}\label{in}
	u(0,x)= u_0(x) \in  (L^1 \cap L^{\infty})(\R^d).
\end{equation}

In \eqref{eq1},
\begin{equation}\label{fregul}
	f= (f_1,\ldots, f_d) \in W^{1,\infty}(\R;\R^d)
\end{equation}
is a given vector-valued flux function, $a =(a_{ij})\ge 0$ is a
given symmetric matrix-valued diffusion function of the form
\begin{equation}\label{dif}
	a= \sigma^a (\sigma^a)^{\mathrm{tr}}, \qquad
	\sigma^a \in \R^{d \times K},
	\quad 1 \leq K \leq d.
\end{equation}
More precisely, the components of $a$ are
$a_{ij} = \sum_{k=1}^K \sigma_{ik}^a \sigma_{jk}^a$ for $i,j= 1,\ldots,d$.
We assume that the matrix-valued function
$\sigma^a = (\sigma_{ik}^a):\R \to \R^{d \times K}$ satisfies
\begin{equation}\label{sigreg}
	\sigma^a\in W^{1,\infty}(\R;\R^{d\times K}).
\end{equation}
Observe that we do not assume the matrix $a(\cdot)$ to be 
strictly positive definite, so the operator
$\Div (a(u) \nabla u)$ may be strongly degenerate, and
hence the phrase ``mixed hyperbolic-parabolic" is justified.

In terms of its singular integral representation, the nonlocal
operator $\Levy$ in \eqref{eq1} takes the form
\begin{equation} \label{nlocsplit}
	\Levy[u](t,x) =  \int_{\R^d \setminus \{0\}}
	\left[ u(t,x+ z) - u(t,x) - z \cdot \nabla u\, \mathbf{1}_{|z|<1} 
	\right]\, \pi(dz),
\end{equation}
where the singular L\'evy measure $\pi(dz)$ is a positive, $\sigma$-finite Borel measure
on $\R^d \setminus \{0\}$ satisfying $\pi(\{0\}) = 0$, $\pi(d(-z))=-\pi(dz)$, and
\begin{equation}\label{intcond1}
	\int_{\R^d \setminus \{0\}} \left( \abs{z}^2 \mathbf{1}_{\abs{z}<1}
	+ \abs{z} \mathbf{1}_{\abs{z}\ge 1}\right)\, \pi(dz) < \infty,
\end{equation}
where we note that $z$ can be replaced by a certain 
regular jump function $j(z)$ easily throughout the analysis.
A typical example is provided by taking
$$
\pi(z)=\frac{1}{|z|^{d+\alpha}}\, 
\mathbf{1}_{|z|<1} \dz,\qquad \alpha \in (0,2).
$$
This example corresponds to the fractional
Laplacian $\Delta_\alpha := -(-\Delta)^{\frac{\alpha}{2}}$
on $\R^d$, which can also be defined in terms of the Fourier transform as
$$
\widehat{\Delta_\alpha v}(\omega)= \abs{\omega}^{\alpha}
\widehat{v}(\omega), \qquad \omega\in \R^d.
$$
This definition is employed in \cite{Droniou:2006os} to prove  (\ref{nlocsplit}) in this case.

Nonlocal operators like $\Delta_\alpha$ are examples of a
pseudodifferential operator $\mathcal{P}$ with a symbol $a(\omega)\ge 0$ such that
$\widehat{\mathcal{P}v}(\omega)=a(\omega) \widehat{v}(\omega)$.
The function $e^{-ta(\omega)}$ is positive definite, and thus, by the
L\'evy-Khintchine formula, it can be represented as
$$
a(\omega)=i b\cdot \omega + q(\omega)
+  \int_{\R^d \setminus \{0\}} \left(1-e^{-i z\cdot\omega}
- i z\cdot \omega\, \mathbf{1}_{\abs{z}<1}(z)\right)\,\pi(dz),
$$
where $b\in \R^d$ represents the drift term, $q(\omega)=\sum_{i,j=1}^d q_{ij} \omega_i \omega_j$
is a positive definite quadratic function representing the pure diffusion
part ($q(\omega)=\abs{\omega}^2$ gives raise to the usual Laplacian $-\Delta$),
and the L\'evy measure $\pi(dz)$ accounts for the jump (non-local) part.
In our setting of $\Levy$, cf.~\eqref{nlocsplit}, we assume
$b\equiv 0$  and $q\equiv 0$, i.e, we are dealing with a pure jump
operator. The key point is that any pseudodifferential operator $\mathcal{P}$ is the generator
of a L\'evy process which is completely characterized in terms of the triplet $\left(b,q,\pi(dz)\right)$.
For more details about the L\'evy-Khintchine formula and L\'evy processes in general, we
refer to \cite{Bertoin:1996ud,Jacob:2001rf,Jacob:2002xp,Jacob:2005ss,Sato:1999qd}.

Integro-partial differential equations, also known as nonlocal, fractional
or L\'evy partial differential equations, appear frequently in many
different areas of research and find many applications in engineering and finance,
including nonlinear acoustics, statistical mechanics, biology, fluid flow,
pricing of financial instruments, and portfolio optimization.
Many authors have recently contributed to advancing the mathematical theory for quasilinear and
fully nonlinear partial differential equations that are supplemented with a
fractional diffusion operator arising as the generator of a L\'evy semigroup, addressing
questions like existence, uniqueness, regularity, formation of
singularities, and asymptotic behavior of solutions.

For results with reference to fully nonlinear equations, such as the Hamilton-Jacobi-Bellman equation,
and the (in this context relevant) theory of viscosity solutions, we refer to \cite{Alibaud:2007dw,Alvarez:1996lq,Arisawa:2006db,Barles:1997xj,Barles:2008lq,Caffarelli:2007hh,Caffarelli:2007qr,Caffarelli:2008bx,Garroni:2002il,Jakobsen:2005jy,Jakobsen:2006aa,Pham:1998zt,Sayah:1991jk,Sayah:1991xy,Silvestre:2006bq,Silvestre:2007qp,Soner:IntegroDif}, see also
\cite{Benth:2001tv,Benth:2002mk,Cont:2004gk} for some concrete applications to finance.

More recently, a number of authors \cite{Alibaud:2007mi,Alibaud:2007qe,Biler:1998ai,Biler:2001th,Bossy:2002eu,Brandolese:2008sp,Droniou:2006os,Karch:2008dp} have studied questions regarding
existence, uniqueness, regularity, and temporal asymptotics for
quasilinear equations, such as the fractal Burgers equation
\begin{equation}\label{eq:fractional-Burgers}
	\pt u + \partial_x (u^2/2)=-(-\partial_{xx}^2)^{\frac{\alpha}{2}}u,
\end{equation}
and more generally multi-dimensional fractional conservation laws
\begin{equation}\label{eq:fCL}
	\partial_t u + \Div f(u) =  \Delta_\alpha u,
\end{equation}
where the parameter $\alpha$ is assumed lie in the interval $(0,2)$.
Of course,  the excluded case $\alpha=2$ corresponds to the already fully understood viscous
conservation law  $\pt u + \Div f(u)=\Delta u$, solutions of which
are always smooth in $t>0$. Regarding the less studied case $\alpha\in [1,2)$, it was recently proved
in \cite{Droniou:2003mz,Kiselev:2008jt} that solutions of the fractional
Burgers equation \eqref{eq:fractional-Burgers} are also smooth in $t>0$.
In the case $\alpha<1$ for the fractional conservation
law \eqref{eq:fCL} the order of the diffusion part
is lower than the first order hyperbolic part, so we do not expect any regularizing
effect to take place. Indeed, for the fractional Burgers
equation \eqref{eq:fractional-Burgers} with $\alpha<1$ it is proved
in \cite{Alibaud:2007qe,Kiselev:2008jt} that solutions can
develop discontinuities in finite time. Consequently, one should employ
a notion of entropy solutions for fractional conservation
laws \eqref{eq:fCL}, i.e., weak solutions satisfying an
additional entropy condition, to ensure the global-in-time well-posedness.
This is well-known for conservation
laws $\pt u + \Div f(u)=0$, cf.~Kru\u{z}kov \cite{Kruzkov:1970kx}, and
the well-posedness theory of Kru\u{z}kov was recently extended to
fractional conservation laws in \cite{Alibaud:2007mi}.

In recent years the theory of Kru\u{z}kov \cite{Kruzkov:1970kx} has been extended to
quasilinear mixed hyperbolic-parabolic equations of the form
\begin{equation}\label{eq1-local}
	\partial_t u + \Div f(u) = \Div (a(u) \nabla u),
\end{equation}
where $f$ and $a$ satisfy \eqref{fregul} and \eqref{dif}-\eqref{sigreg}, respectively. Since
the diffusion matrix $a(u)$ is not assumed to be strictly positive definite, \eqref{eq1-local}
is strongly degenerate and will in general posses discontinuous solutions.
In the isotropic case (with $a(\cdot)$ being a scalar function) the
first general uniqueness result is due to Carrillo \cite{Carrillo:1999hq}, who developed
an original extension of Kru\u{z}kov's method of doubling variables to prove
his result, cf.~\cite{Karlsen:2002bh,Karlsen:2003za,Mascia:2002dq,Michel:2003tw}
for some additional applications of his techniques. The anisotropic
case ($a(\cdot)$ being a matrix-valued function) was first treated by Chen
and Perthame \cite{Chen:2003td}, who developed a
kinetic formulation and established the uniqueness result using regularization by convolution.
An alternative proof of the result of Chen and Perthame, adapting
the device of doubling variables, was developed
in \cite{Bendahmane:2004go}, cf.~also \cite{Chen:2006oy,Chen:2005wf,Perthame:2003yq}
some other papers dealing with the anisotropic case.

The main purpose of this paper is to extend the uniqueness and
``continuous dependence on the nonlinearities" results of
\cite{Bendahmane:2004go,Chen:2006oy,Chen:2005wf,Perthame:2003yq} to
fractional degenerate parabolic equations of the form \eqref{eq1}.
We introduce the notion of entropy solutions and state
the main results in Section \ref{sec:entform} .
Sections \ref{prexis} (existence), \ref{pruniq} (uniqueness),
and \ref{prcdep} (continuous dependence on the nonlinearities
and the L\'evy measure) are devoted to the proofs of the main results.

\section{Notion of solution and main results}\label{sec:entform}

For $i =1,\ldots,d$ and $k=1,\ldots,K,$ define
$$
\zeta_{ik}^a(z) := \int_0^z \sigma_{ik}^a(\xi)\dxi,
\quad \zeta_{ik}^{a,\psi}(z) = 
\int_0^z \psi(\xi) \sigma_{ik}^a(\xi)\dxi,
\quad  z\in \R,
$$
for any $\psi \in C(\R)$. Given any convex $C^2$
entropy function $\eta:\R \to \R$, we define
the corresponding entropy fluxes $q=(q_i):\R \to \R^d$
and $r=(r_{ij}):\R \to \R^{d \times d}$ by
$$
q'(z) = \eta'(z) f(z), \qquad r'(z) = \eta'(z) a(z).
$$
We refer to $(\eta, q, r)$ as an entropy-entropy flux triple.

We now introduce the entropy formulation of  \eqref{eq1}-\eqref{in}.

\begin{definition}\label{def:entropy}
An entropy solution of the initial value problem \eqref{eq1}-\eqref{in}
is a measurable function $u : Q_T \to \R$ satisfying the following conditions:

\begin{Definitions}
\item $u \in L^{\infty}(Q_T)$, $u \in L^{\infty}(0,T;L^1(\R^d))$,
\begin{equation}\label{l2regu}
\sum_{i=1}^d \pxi \zeta_{ik}^a(u)
\in L^2(Q_T), \qquad k=1,\ldots,K,
\end{equation}
and
\begin{equation}\label{frintcnd}
	\iint_{Q_T}\int_{\R^d \setminus \{0\}}
	\left(u(t,x+z) - u(t,x) \right)^2 \,
	\pi (dz) \dx \dt < +\infty.
\end{equation}

\item For $k= 1,\ldots, K$,
\begin{equation}\label{crule}
	\sum_{i=1}^d \pxi \zeta_{ik}^{a,\psi}(u)
	= \psi(u) \sum_{i=1}^d \pxi \zeta_{ik}^a(u),
	\quad \text{a.e.~in $Q_T$ and in $L^2(Q_T)$},
\end{equation}
for any $\psi \in C(\R)$.
\label{def:chainrule}

\item For any entropy-entropy flux triple $(\eta, q, r)$,
\begin{equation}\label{entineq}
	\begin{split}
		& \iint_{Q_T} \Bigl ( \eta(u) \pt \varphi
		+ \sum_{i=1}^d q_i(u) \pxi \varphi
		+ \sum_{i,j=1}^d r_{ij}(u) \pxixj \varphi \Bigr) \dx\dt
		\\ & \qquad
		+ \iint_{Q_T} \eta(u) \Levy[\varphi] \dx\dt
		+ \int_{\R^d} \eta(u_0) \varphi(0,x) \dx
		\geq  n^u + m^u,
	\end{split}
\end{equation}
for all non-negative $\varphi \in C_c^{\infty}([0,T) \times \R^d)$, where
\begin{align*}
	n^u & =  \iint_{Q_T} \eta''(u) \sum_{k=1}^K
	\Bigl(\sum_{i=1}^d \pxi \zeta_{ik}^a(u) \Bigr)^2 \varphi(t,x)\dx \dt,  \\
	m^u & = \iint_{Q_T}  \int_{\R^d \setminus \{0\}}
	\overline{\eta''}(u;z)
	\left(u(t,x+z)-u(t,x)\right)^2 \varphi(t,x)  \, \pi(dz)\dx\dt,
\end{align*}
and
$$
\overline{\eta''}(u;z)
= \int_0^1 (1-\tau) \eta''((1-\tau)u(t,x)+\tau u(t,x+z)) \dtau.
$$
If, in addition,
\begin{equation}\label{eq:add-cond}
	\left\{
	\begin{split}
		& |z|\, \pi(dz) \in L^1(\R^d \setminus \{0\}),\\
		\text{or}\\
		& |z|^2\, \pi(dz) \in L^1(\R^d \setminus \{0\}) 
		\quad \text{and} \quad u\in L^\infty(0,T;BV(\R^d)),
	\end{split}\right.
\end{equation}
then we can drop the fractional parabolic dissipation measure $m^u$ and 
replace \eqref{entineq} by the simpler condition
\begin{equation}\label{entineq-simpler}
	\begin{split}
		& \iint_{Q_T} \Bigl ( \eta(u) \pt \varphi
		+ \sum_{i=1}^d q_i(u) \pxi \varphi
		+ \sum_{i,j=1}^d r_{ij}(u) \pxixj \varphi \Bigr) \dx\dt
		\\ & \qquad
		+ \iint_{Q_T} \eta(u) \Levy[\varphi] \dx\dt
		+ \int_{\R^d} \eta(u_0) \varphi(0,x) \dx\geq  n^u.
	\end{split}
\end{equation}
\end{Definitions}
\end{definition}

We remark that the chain rule \eqref{crule}
is automatically fulfilled when $a(\cdot)$ is a scalar or
a diagonal matrix, cf.~Chen and Perthame \cite{Chen:2003td}, and in
this case we can drop \ref{def:chainrule} from the definition.

Starting from the definition of $\Levy$ (cf.~calculations in 
the upcoming sections),  we can replace the term
$$
\iint_{Q_T} \eta(u) \Levy[\varphi] \dx\dt-m^u,
$$
occurring in \eqref{entineq} by 
\begin{align*}
	& \iint_{Q_T} \int_{|z|< r}\eta(u) [\varphi(t,x+z)-\varphi(t,x)
	-\nabla \varphi \cdot z]\, \pi(dz) \dx\dt,
	\\ & \quad 
	+ \iint_{Q_T} \int_{|z|\geq r}\eta'(u) [u(t,x+z) - u(t,x)]\, \pi(dz) \dx\dt, 
	\quad \forall r\in (0,1),
\end{align*}
This formulation of the nonlocal term is directly related 
to the formulation used in \cite{Alibaud:2007mi} 
for fractional conservation laws.

Our first result is the expected $L^1$ contraction
property (and thus the uniqueness) of entropy solutions.

\begin{theorem}\label{th:unique}
Suppose $f$ and $a$ satisfy \eqref{fregul} 
and \eqref{dif}-\eqref{sigreg}, respectively, and
that the L\'evy measure $\pi(dz)$
satisfies \eqref{intcond1}. Then there exists an
entropy solution of \eqref{eq1}-\eqref{in}.
Let $u,v$ be two entropy solutions of \eqref{eq1} with
initial data $u|_{t=0} =u_0 \in (L^1\cap L^{\infty})(\R^d)$,
$v|_{t=0} = v_0 \in (L^1\cap L^{\infty})(\R^d)$.
For a.e. $t \in (0,T),$  we have
\begin{equation}\label{eq:contraction}
	\int_{\R^d} \left(u(t,x)-v(t,x)  \right)^+ \dx
	\leq \int_{\R^d} \left(u_0-v_0 \right)^+ \dx.
\end{equation}
Consequently, if $u_0 \leq v_0$ a.e.~in $\R^d$ then $u \leq v$ a.e.~in $Q_T$, so
whenever $u_0=v_0$ a.e.~in $\R^d$, then $u=v$ a.e.~in $Q_T$. 

Finally, if we replace \eqref{entineq} by \eqref{entineq-simpler}, then the 
$L^1$ contraction principle \eqref{eq:contraction} 
continues to hold provided \eqref{eq:add-cond} is effective; it 
is sufficient that (say) only $v$ belongs to $L^\infty(0,T;BV(\R^d))$ in the case 
$\int |z|^2 \,\pi(dz)<\infty$.
\end{theorem}

This theorem generalizes to the ``non-local diffusion" case
the result of Chen and Perthame \cite{Chen:2003td}. The proof follows that
of Bendahmane and Karlsen \cite{Bendahmane:2004go}. 

Regarding the last part of Theorem \ref{th:unique}, assuming $v_0\in BV(\R^d)$ it follows 
from \eqref{eq:contraction} that $v\in L^\infty(0,T;BV(\R^d))$, as required.

Our second result,  which is a refinement of the previous theorem, reveals
how the entropy solution $u$ depends on  the L\'evy measure $\pi(dz)$, and the
nonlinear fluxes $f,a$ (i.e., it is a ``continuous dependence" estimate).

\begin{theorem}\label{th:contdep}
Suppose $f$ and $a$ satisfy \eqref{fregul} and 
\eqref{dif}-\eqref{sigreg}, respectively, and
that the L\'evy measure $\pi(dz)$ satisfies \eqref{intcond1}.
Let $u \in L^{\infty}(0,T;BV(\R^d))$ be the entropy 
solution of \eqref{eq1} with $BV$ initial 
data $u_0 \in (L^1\cap L^{\infty} \cap BV)(\R^d)$ 
and with a L\'evy measure of the form $\pi(dz) = m(z)\, dz$ 
for some integrable function $m:\R^d \setminus \{0\}\to \R_+$.

Replace the data set
$$
(f,a,\pi,u_0), \quad a=\sigma^a(\sigma)^{\mathrm{tr}}, 
\quad \pi(dz)=m(z)\,dz
$$
by another data set
$$
(\tilde f,\tilde a,\tilde \pi(dz),v_0), \quad
\tilde a=\sigma^{\tilde a}(\sigma^{\tilde a})^{\mathrm{tr}}, 
\quad \tilde \pi(dz)=\tilde m(z)\,dz,
$$
where $\tilde f,\sigma^{\tilde a}, \tilde\pi, \tilde m$ satisfy the same regularity conditions
as $f, \sigma^a, \pi, m$ and moreover $v_0\in (L^1\cap L^\infty)(\R^d)$.
Denote the corresponding entropy solution by $v$, and assume that $v\in C([0,T];L^1(\R^d))$. 
Suppose $u$ and $v$  take values in a closed interval $I \subset \R$.

For any $t \in (0,T)$,
\begin{equation}\label{contdest}
	\begin{split}
		&\norm{u(t,\cdot) - v(t,\cdot)}_{L^1(\R^d)}
		\\ &\leq \norm{u_0-v_0}_{L^1(\R^d)}
		+ C_1 t \norm{f-\tilde f}_{W^{1,\infty}(I);\R^d)}
		+ C_2 \sqrt{t} \norm{\sigma^a-\sigma^{\tilde a}}_{L^{\infty}(I;\R^{d \times K})} \\
		& \quad
		+ C_3 \sqrt{t}\sqrt{\left(\int_{|z|<1}\abs{z}^2 \abs{m(z)-\tilde m(z)} \dz\right)}
		+ C_4 t\int_{|z|\geq 1}|z| \abs{m(z)-\tilde m(z)} \dz,
	\end{split}
\end{equation}
where the constants $C_i$, $i=1,\ldots,4$, depend on 
the $L^{\infty}(0,T;BV(\R^d))$ norm of $u$.
\end{theorem}

This theorem generalizes results in \cite{Chen:2005wf,Chen:2006oy} 
to the ``fractional case". 

\section{Proof of Theorem \ref{th:unique} (existence)}\label{prexis}

Although a detailed version of the existence of entropy 
solutions to (\ref{eq1}) is presented in \cite{KU2}, to 
motivate the entropy condition and to present a brief sketch, let 
us consider the following accompanying problem
containing a uniformly parabolic operator depending
on a small parameter $\rho>0$:
\begin{equation}\label{eq1-visc}
	\partial_t u_\rho + \Div f(u_\rho) =
	\Div( a(u_\rho) \nabla u_\rho)+  \Levy[u_\rho(t,\cdot)]
	+\rho \Delta u_\rho.
\end{equation}
It is standard to construct a smooth solution $u_\rho$
to \eqref{eq1-visc}, for each fixed $\rho>0$. 
Indeed, it can be done using the Galerkin method
and the compactness argument, see Chapter 5 in \cite{Ev} and \cite{Kiselev:2008jt}.

As usual, the game is to pass
to the limit as $\rho \to 0$ and identify the entropy
condition satisfied by the limit function $u$. 
We will be brief in establishing the following estimates, since 
most of them are similar to the ones in \cite{Chen:2003td} and 
we will assume $u_0 \in W^{2,1}\cap H^1 \cap L^{\infty}(\R^d)$, 
for general $u_0 \in L^1(\R^d)$ one can follow the 
approximation procedure presented in \cite{Chen:2003td}.

The following estimates can be established for sufficiently regular initial data:
$$
\norm{u_\rho}_{L^\infty(Q_T)}\le C;\qquad \abs{u_\rho(t,\cdot)}_{BV(\R^d)}\le C;
$$
$$
\norm{u_\rho(t_2,\cdot)-u_\rho(t_1,\cdot)}_{L^1(\R^d)}\to 0,
\quad \text{as $\abs{t_2-t_1}\to 0$, uniformly in $\rho$.}
$$
Hence there is a limit $u$ such that, passing if necessary to a subsequence as $\rho\to 0$,
\begin{equation}\label{eq:ueps-conv}
	u_\rho\to u\quad \text{a.e.~in $Q_T$ and in $L^p(Q_T)$ for any $p \in [1,\infty)$.}
\end{equation}

Next, we derive an energy estimate. To this end, fix a convex $C^2$ function $\eta$ and
define $q,r$ by $q'=\eta'f'$, $r'=\eta' a$. Multiplying \eqref{eq1-visc} by $\eta'$ yields
\begin{equation}\label{eq1-visc-entropy}
	\partial_t \eta(u_\rho) + \Div q(u_\rho)
	=\sum_{i,j=1}^d \partial_{ij}^2 r_{ij}(u_\rho) + \Levy[\eta(u_\rho]
	+ \rho \Delta \eta(u_\rho) - \nu_\rho
\end{equation}
where $\nu_\rho=\nu_\rho^1+\nu_\rho^2+\nu_\rho^3$ consists of
three parts:

(i) the entropy dissipation term
$$
\nu_\rho^1:=\rho\Delta \eta(u_\rho)-\rho\eta'(u_\rho)\Delta u_\rho
=\rho \eta''(u_\rho) \abs{\nabla u_\rho}^2;
$$

(ii) the parabolic dissipation term
$$
\nu_\rho^2:= \sum_{i,j=1}^d \partial_{ij}^2 
r_{ij}(u_\rho)- \eta'(u_\rho)\Div (a(u_\rho) \nabla u_\rho)
=\eta''(u_\rho) \sum_{k=1}^K 
\left( \sum_{i=1}^d \partial_{x_i} \zeta_{ik}^a(u_\rho)\right)^2;
$$

(iii) the fractional parabolic dissipation term	
$$
\nu_\rho^3 =  \int_{\R^d \setminus \{0\}}\overline{\eta''}(u_\rho;z)
\left(u_\rho(t,x+z)-u_\rho(t,x)\right)^2\, \pi(dz),
$$
where $\overline{\eta''}(u_\rho;z)= \int_0^1 (1-\tau) 
\eta''((1-\tau)u_\rho(t,x) + \tau u_\rho(t,x+z)) \dtau$.

In deriving \eqref{eq1-visc-entropy}, the ``new" computation
is the one showing that the commutator
$$
\Levy[\eta(u_\rho)]-\eta'(u_\rho) \Levy[u_\rho]
$$
equals $\nu_\rho^3$, but this follows easily from
Taylor's formula with integral reminder:
 \begin{equation}\label{eq:Taylor}
 	\begin{split}
		\eta(b)-\eta(a) &= \eta'(a) \left(b - a\right)
		\\ & \qquad
		+\left(\int_0^1 (1-\tau) \eta''((1-\tau)a + \tau b) \dtau\right)
		\left(b - a\right)^2.
	\end{split}
\end{equation}

Specifying $\eta(z)=\frac{z^2}{2}$ in \eqref{eq1-visc-entropy} gives
$$
\int_0^T \int_{\R^d}  \sum_{k=1}^K \left( \sum_{i=1}^d \partial_{x_i}
\zeta_{ik}^a(u_\rho)\right)^2\dx \dt \le C
$$
and
\begin{equation}\label{eq:weakconv-zeta}
	\sum_{i=1}^d \partial_{x_i} \zeta_{ik}^a(u_\rho) 
	\weak \sum_{i=1}^d \partial_{x_i} \zeta_{ik}^a(u) 
	\quad \text{in $L^2(Q_T)$}.
\end{equation}
From this we easily see, as in \cite{Chen:2003td}, that 
\eqref{l2regu} and \eqref{crule} in Definition \ref{def:entropy} hold.

Regarding the non-local operator $\Levy$, the same choice for $\eta$ 
reveals that (\ref{frintcnd}) in Definition \ref{def:entropy} holds. Now set
$$
\Pi(dz):= \left( \abs{z}^2 \mathbf{1}_{|z|<1}
+\abs{z}\, \mathbf{1}_{|z|\ge 1}\right)\, \pi(dz),
$$
and note that $\Pi(dz)$ is a bounded Radon measure.
Introducing the short-hand notation
$$
D_\rho(t,x,z)= \frac{u_\rho(t,x+z) 
- u_\rho(t,x)}{\abs{z}\mathbf{1}_{_{|z|<1}}
+ \sqrt{\abs{z}} \mathbf{1}_{_{|z|\ge 1}}}
\qquad
d\mu=\Pi(dz)\otimes dx\otimes dt,
$$
\eqref{frintcnd} translates into $D_\rho$ being uniformly bounded
in $L^2((0,T)\times \R^d\times (\R^d \setminus \{0\});d\mu)$. Consequently, we may
assume that there is a limit function $D$ such that
$$
D_\rho \weak D \quad \text{in $L^2((0,T)\times \R^d\times (\R^d \setminus \{0\});d\mu)$.}
$$

Let us identify $D$. To this end, fix a smooth function
$\varphi$ in $C^\infty_c(Q_T)$ and observe
\begin{align*}
	&\iint_{Q_T}\int_{\R^d \setminus \{0\}}
	\varphi(t,x)\frac{u_\rho(t,x+z)-u_\rho(t,x)}{\abs{z}\mathbf{1}_{|z|<1}
	+ \sqrt{\abs{z}} \mathbf{1}_{|z|\ge 1}}
	\, \Pi(dz)\, dx\, dt
	\\ & = \iint_{Q_T}\int_{\R^d \setminus \{0\}}
	\frac{\varphi(t,x+z)-\varphi(t,x)}{\abs{z}\mathbf{1}_{|z|<1}
	+ \sqrt{\abs{z}} \mathbf{1}_{|z|\ge 1}} u_\rho(t,x)
	\, \Pi(dz)\, dx\, dt.
\end{align*}
Now, using that $u_\rho \overset{\rho\to 0}{\longrightarrow} u$ a.e.~in $Q_T$, we
conclude that
$$
D_\rho \weak \frac{u(t,x+z) - u(t,x)}{\abs{z}\mathbf{1}_{|z|<1}
+ \sqrt{\abs{z}} \mathbf{1}_{|z|\ge 1}}
\quad \text{in $L^2((0,T)\times \R^d\times 
(\R^d \setminus \{0\});d\mu)$.}
$$

We are now in a position to pass to the distributional 
limit in \eqref{eq1-visc-entropy} to recover the desired
entropy condition satisfied by the limit 
$u=\lim_{\rho \to 0} u_\rho$.  Note that to
interpret \eqref{eq1-visc-entropy} in the sense of distributions
we use the formula
\begin{equation}\label{eq:IBP}
	\int_{\R^d} \Levy[\Phi(x)] \phi(x) \dx = \int_{\R^d} \Phi(x) \Levy[\phi(x)] \dx,
\end{equation}
which holds for all sufficiently regular (say, $C^2$) functions $\Phi,\phi:\R^d\to \R$.
This relation is easily obtained by a change of variables
$(t,x,z) \mapsto (t,x+z,-z)$ and an integration by parts in $x$.

We claim that the entropy condition satisfied by the limit 
$u=\lim_{\rho \to 0} u_\rho$ takes the following form:
for any convex $C^2$ entropy function $\eta$ and corresponding entropy
fluxes $q,r$ defined by $q'=\eta'f', r'=\eta' a$,
\begin{equation}\label{eq1-limit-entropy}
	\partial_t \eta(u) + \Div q(u)
	\le \sum_{i,j} \partial_{x_i x_j} r_{ij}(u) + \Levy[\eta(u)] - n^{u,\eta}-m^{u,\eta}
\end{equation}
in the sense of distributions, where
$$
n^{u,\eta} = \eta''(u) \sum_{k=1}^K 
\left(\sum_{i=1}^d \partial_{x_i} \zeta_{ik}^a(u)\right)^2
$$
is the parabolic dissipation measure with respect to $u$ and 	
$$
m^{u,\eta} = \int_{\R^d \setminus \{0\}}\overline{\eta''}(u;z)
\left(u(t,x+z) - u(t,x) \right)^2\, \pi(dz),
$$
is the fractional parabolic dissipation measure with respect to $u$.

In view of \eqref{eq:ueps-conv}, to verify \eqref{eq1-limit-entropy}
we only need to argue that
$$
\liminf_{\rho \to 0} \iint_{Q_T} \nu_\rho \, dx \, dt \ge 
\iint_{Q_T} \left( n^{u,\eta} + m^{u,\eta} \right) \dx \dt.
$$
First, $\iint_{Q_T} \nu_\rho^1 \dx \dt \ge 0$ for each $\rho>0$. Second, thanks 
to the weak convergence \eqref{eq:weakconv-zeta} and a standard
weak lower semi-continuity result for quadratic functionals,
\begin{align*}
	& \liminf_{\rho\to 0}
	\int_0^T \int_{\R^d} \eta''(u_\rho)\sum_{k=1}^K 
	\left( \sum_{i=1}^d \partial_{x_i} \zeta_{ik}^a(u_\rho)  \right)^2 
	\varphi \dx \dt \\ & \qquad
	\ge \int_0^T \int_{\R^d} \eta''(u)\sum_{k=1}^K \left( \sum_{i=1}^d 
	\partial_{x_i} \zeta_{ik}^a(u)  \right)^2 \varphi \dx \dt,
\end{align*}
for all test functions $\varphi\in C^\infty_c$. Similarly,
\begin{align*}
	&\liminf_{\rho\to 0}
	\iint_{Q_T}  \int_{\R^d \setminus \{0\}}
	\overline{\eta''}(u_\rho;z)
	\left(u_\rho(t,x+z)-u_\rho(t,x)\right)^2\varphi\, \pi(dz)\dx\dt
	\\ & \qquad
	\ge \iint_{Q_T} \int_{\R^d \setminus \{0\}}\overline{\eta''}(u;z)
	\left(u(t,x+z)-u(t,x)\right)^2\varphi\, \pi(dz)\dx\dt,
\end{align*}
for all test functions $\varphi\in C^\infty_c$. 
Combining, we deduce that \eqref{entineq} 
in Definition \ref{def:entropy}  holds. This completes the proof.

\section{Proof of Theorem \ref{th:unique} (uniqueness)}\label{pruniq}

We shall need $C^2$ approximations $\eta_{\eps}^\pm(z)$ of the functions
$$
\eta^\pm(z) :=(z)^\pm=\max\left(\pm(z),0\right),
\qquad z\in \R.
$$
We build these by picking nondecreasing $C^1$
approximations $\sgn_{\eps}^\pm(z)$ of
$$
\sgn^+(z):=
\begin{cases}
	0, & \textrm{if $z \leq 0$},\\
	1, & \textrm{if $z>0$},
\end{cases}
\qquad
\sgn^-(z) :=
\begin{cases}
	-1, & \textrm{if $z \leq 0$,}\\
	0, & \textrm{if $z>0$,}
\end{cases}
$$
and defining
$$
\eta_{\eps}^\pm(z)
:= \int_0^z \sgn_{\eps}^\pm(\xi)  \dxi,
\qquad z\in \R.
$$
For example, we can take
$$
\sgn_{\eps}^+(z)=
\begin{cases}
	0, & \textrm{if $z<0$},\\
	\sin(\frac{\pi}{2\eps}z), & \textrm{if $0\le z \leq \eps$},\\
	1, & \textrm{if $z>\eps$}.
\end{cases}
\quad
\sgn_{\eps}^-(z)=
\begin{cases}
	-1, & \textrm{if $z<-\eps$},\\
	\sin(\frac{\pi}{2\eps}z), & \textrm{if $-\eps\le  z \le 0$},\\
	0, & \textrm{if $z>0$}.
\end{cases}
$$

The functions $\eta_{\eps}^\pm$
are $C^2$ and convex. Moreover,
$$
\eta_{\eps}^\pm(z)
\overset{\eps\to 0}{\longrightarrow}
\eta^\pm(z),
\qquad z\in \R.
$$

Observe that $\left(\eta_{\eps}^\pm(\cdot-c)\right)_{c \in\R}$
is a family of entropies.
Given these entropies, we introduce the
corresponding entropy fluxes
\begin{align*}
	q_{\eps}^\pm(z,c) & = \int_c^z (\eta_\eps^\pm)'(\xi-c) f'(\xi) d\xi, \qquad z,c \in \R,\\
	r_{\eps}^\pm(z,c)  & = \int_c^z (\eta_{\eps}^\pm)'(\xi-c) a(\xi) \dxi,  \qquad z,c \in \R.
\end{align*}
Clearly, as $\eps \to 0$,
\begin{align*}
	q_{\eps}^\pm(z,c) \to q^\pm(z,c) & := \sgn^\pm(z-c) (f(z)-f(c)), \qquad z,c\in \R,\\
	r_{\eps}^\pm(z,c) \to r^\pm(z,c) & := \sgn^\pm(u-c) (A(u) - A(c)), \qquad z,c\in \R,
\end{align*}
where the (matrix-valued) function $A(\cdot)$
is defined by $\dsp A(z) = \int_0^u  a(\xi) \dxi$.

Observe that $\left(\eta_{\eps}^\pm(\cdot-c),q_{\eps}^\pm(\cdot,c),r_{\eps}^\pm(\cdot,c)\right)_{c \in \R}$
is a family of entropy-entropy flux triples, so choosing
$\eta = \eta_{\eps}^\pm$ in \eqref{entineq} yields
\begin{equation}\label{in1}
	\begin{split}
		& \iint_{Q_T} \Bigl( \eta_{\eps}^\pm(u-c) \pt \varphi
		+  \sum_{i=1}^d  q_{\eps,i}^\pm(u,c) \pxi \varphi
		+ \sum_{i,j=1}^d  r_{\eps,ij}^\pm(u,c) \pxixj \varphi \Bigr)\dx \dt
		\\ & \quad + \iint_{Q_T}  \eta_{\eps}^\pm(u-c) \Levy[\varphi]\dx \dt
		+ \int_{\R^d} \eta^\pm(u_0-c) \varphi(0,x) \dx
		\\ & \geq  \iint_{Q_T} (\eta_{\eps}^\pm)''(u-c)
		\sum_{k=1}^K \Bigl(\sum_{i=1}^d  \pxi \zeta_{ik}^a(u)\Bigr)^2 \varphi \dx \dt \\
		&\quad + \iint_{Q_T} \int_{\R^d \setminus \{0\}}
		\overline{(\eta_{\eps}^\pm)''}(u-c;z)
		\left(u(t,x+z)-u(t,x)\right)^2 \varphi \, \pi(dz) \dx \dt.
	\end{split}
\end{equation}
Moreover,
\begin{align*}
	\overline{(\eta_\eps^\pm)''}(u-c;z)
	& = \int_0^1 (1-\tau) (\eta_\eps^\pm)''\Bigl((1-\tau)u(t,x)+\tau u(t,x+z),c\Bigr) \dtau
	\\ & =  \int_0^1 (1-\tau) (\sgn_\eps^\pm)'\Bigl((1-\tau)(u(t,x)-c)+\tau (u(t,x+z)-c)\Bigr)\dtau.
\end{align*}

To proceed, the following simple observations will be useful:
\begin{itemize}
	\item  $\sgn_{\eps}^-(z-c) =-\sgn_{\eps}^+(c-z)$ and $\eta_{\eps}^-(z-c) = \eta_{\eps}^+(c-z)$;
	\item  $q_{\eps}^-(z,c) = q_{\eps}^+(c,z)$ and $r_\eps^-(z,c) = r_\eps^+(c,z)$;
	\item  $(\eta_{\eps}^-)''(z-c) = (\eta_{\eps}^+)''(c-z)$.
\end{itemize}
Employing these observations, we can rewrite the ``$-$" part of \eqref {in1} as
\begin{equation}\label{in3}
	\begin{split}
		& \iint_{Q_T} \Bigl( \eta_{\eps}^+(c-u) \pt \varphi
		+  \sum_{i=1}^d  q_{\eps,i}^+(c,u) \pxi \varphi
		+ \sum_{i,j=1}^d  r_{\eps,ij}^+(c,u) \pxixj \varphi \Bigr)\dx \dt
		\\ & \quad + \iint_{Q_T}  \eta_{\eps}^+(c-u) \Levy[\varphi]\dx \dt
		+ \int_{\R^d} \eta_{\eps}^+(c-u_0) \varphi(0,x) \dx
		\\ & \geq  \iint_{Q_T} (\eta_{\eps}^+)''(c-u)
		\sum_{k=1}^K \Bigl(\sum_{i=1}^d  \pxi \zeta_{ik}^a(u)\Bigr)^2 \varphi \dx \dt \\
		&\quad + \iint_{Q_T} \int_{\R^d \setminus \{0\}}
		\overline{(\eta_{\eps}^+)''}(c-u;z)
		\left(u(t,x+z)-u(t,x)\right)^2 \varphi \, \pi(dz) \dx \dt.
	\end{split}
\end{equation}

To establish the $L^1$ contraction property \eqref{eq:contraction} we shall
employ the doubling-of-variables device of Kru\u{z}kov \cite{Kruzkov:1970kx}.
Let $u=u(t,x)$, $v=v(s,y)$ be two entropy solutions as stated
in Theorem \ref{th:unique}. Moreover, let $\varphi=\varphi(t,x,s,y)$ be a
test function in the doubled variables $(t,x,s,y)$.
To simplify the presentation, we introduce the following notation (with $\nabla_{x+y}$ being
short-hand for $\nabla_x+\nabla_y$)
\begin{align*}
	\Levy_x[\varphi] & := \int_{\R^d \setminus \{0\}}
	\left[ \varphi(t,x+z,s,y) - \varphi
	- z \cdot \nabla_x \varphi
	\mathbf{1}_{|z|<1} \right]\, \pi(dz), \\
	\Levy_y[\varphi]  & =  \int_{\R^d \setminus \{0\}}
	\left[ \varphi(t,x,s,y+ z) - \varphi
	-z \cdot \nabla_y \varphi
	\mathbf{1}_{|z|<1} \right]\, \pi(dz), \\
	\Levy_{x+y}[\varphi]  & =  \int_{\R^d \setminus \{0\}}
	\Bigl[ \varphi(t,x+z,s,y+ z)-\varphi
	-z \cdot \nabla_{x+y}\varphi
	\mathbf{1}_{|z|<1} \Bigr]\, \pi(dz),
\end{align*}

In the ``$+$" part of \eqref{in1} written the entropy solution $u(t,x)$ we choose $c=v(s,y)$
and integrate the result over $(s,y)$, obtaining
\begin{equation}\label{ineqt04}
	\begin{split}
		& \fourint  \Bigl( \eta_{\eps}^+(u-v) \pt \varphi
		+  \sum_{i=1}^d  q_{\eps,i}^+(u,v) \pxi \varphi
		+ \sum_{i,j=1}^d  r_{\eps,ij}^+(u,c) \pxixj \varphi \Bigr)\dx \dt \dy\ds
		\\ & \quad + \fourint \eta_{\eps}^+(u-v) \Levy_x[\varphi]\dx\dt\dy\ds
		+ \threeint \eta_{\eps}^+(u_0-v) \varphi(0,x,s,y) \dx\dy\ds
		\\ & \geq  \iint_{Q_T} (\eta_{\eps}^+)''(u-v)
		\sum_{k=1}^K \Bigl(\sum_{i=1}^d
		\pxi \zeta_{ik}^a(u)\Bigr)^2 \varphi \dx\dt\dy\ds \\
		&\quad+ \fourint\!\!\int_{\R^d \setminus \{0\}}
		\overline{(\eta_{\eps}^+)''}(u(t,\cdot)-v;z)
		\left(u(t,x+z)-u(t,x)\right)^2 \varphi \, \pi(dz) \dx\dt\dy\ds.
	\end{split}
\end{equation}

Similarly, in \eqref{in3} written for the entropy solution $v(s,y)$ we choose $c=u(t,x)$ and
integrate over $(t,x)$, thereby obtaining
\begin{equation}\label{ineek05}
	\begin{split}
		& \fourint \Bigl( \eta_{\eps}^+(u-v) \ps \varphi
		+  \sum_{i=1}^d  q_{\eps,i}^+(u,v) \pyi \varphi
		+ \sum_{i,j=1}^d  r_{\eps,ij}^+(u,v) \pyiyj \varphi \Bigr)\dx\dt\dy\ds
		\\ & \; + \fourint  \eta_{\eps}^+(u-v) \Levy_y[\varphi]\dx\dt\dy\ds
		+ \threeint \eta_{\eps}^+(u-v_0) \varphi(t,x,0,y) \dx\dt\dy
		\\ & \geq  \fourint (\eta_{\eps}^+)''(u-v)
		\sum_{k=1}^K \Bigl(\sum_{i=1}^d\pyi \zeta_{ik}^a(v)\Bigr)^2
		\varphi \dx \dt\dy\ds \\
		&\; + \fourint\!\!\int_{\R^d \setminus \{0\}}
		\overline{(\eta_{\eps}^+)''}(u-v(s,\cdot);z)
		\left(v(s,y+z)-v(s,y)\right)^2 \varphi \, \pi(dz) \dx\dt\dy\ds.
	\end{split}
\end{equation}

Adding \eqref{ineqt04} and \eqref{ineek05} yields
\begin{equation}\label{in6}
	I_{\mathrm{time}}(\eps)+I_{\mathrm{conv}}(\eps)+I_{\mathrm{diff}}(\eps)
	+I_{\mathrm{fdiff}}(\eps)+	I_{\mathrm{init}}(\eps)
	\ge I_{\mathrm{diss}}(\eps)+I_{\mathrm{fdiss}}(\eps),
\end{equation}
where
\begin{align*}
	I_{\mathrm{time}}(\eps) & = \fourint \eta_{\eps}^+(u-v) (\pt +\ps) \varphi  \dx\dt\dy\ds \\
	I_{\mathrm{conv}}(\eps) & = \fourint \sum_{i=1}^d  q_{\eps,i}^+(u,v)
	(\pxi + \pyi) \varphi \dx\dt\dy\ds \\
	I_{\mathrm{diff}}(\eps) & = \fourint \sum_{i,j=1}^d  r_{\eps,ij}^+(u,v)
	(\pxixj+\pyiyj) \varphi \dx\dt\dy\ds \\
	I_{\mathrm{fdiff}}(\eps) & = \fourint  \eta_{\eps}^+(u-v)
	\Bigl(\Levy_x[\varphi]+\Levy_y[\varphi]\Bigr)\dx\dt\dy\ds\\
	I_{\mathrm{init}}(\eps) & = \threeint \eta_{\eps}^+(u_0-v) \varphi(0,x,s,y) \dx\dy\ds
	\\ & \qquad\quad
	+\threeint \eta_{\eps}^+(u-v_0) \varphi(t,x,0,y) \dx\dt\dy\\
	I_{\mathrm{diss}}(\eps) & = \fourint  (\eta_{\eps}^+)''(u-v)
	\\ & \qquad \qquad
	\times \sum_{k=1}^K \left[ \Bigl(\sum_{i}^d \pxi \zeta_{ik}^a(u)\Bigr)^2
	+\Bigl(\sum_{i=1}^d\pyi \zeta_{ik}^a(v)\Bigr)^2\right]
	\varphi \dx \dt\dy\ds \\
	I_{\mathrm{fdiss}}(\eps) & =  \fourint\!\!\int_{\R^d \setminus \{0\}}
	\Biggl [ \overline{(\eta_{\eps}^+)''}(u(t,\cdot)-v;z) \left(u(t,x+z)-u(t,x)\right)^2
	\\ & \qquad\qquad\qquad\qquad\qquad
	+ \overline{(\eta_{\eps}^+)''}(u,v(s,\cdot);z) \left(v(s,y+z)-v(s,y)\right)^2 \Biggr]
	\\ & \qquad\qquad\qquad\qquad\qquad\qquad\qquad\qquad\qquad
	\times \varphi \, \pi(dz) \dx\dt\dy\ds.
\end{align*}

In view of the inequality ``$a^2+b^2\ge 2ab$", we have
$I_{\mathrm{diss}}(\eps)  \ge \widetilde I_{\mathrm{diss}}(\eps)$, with
$$
\widetilde I_{\mathrm{diss}}(\eps)
= 2\fourint  (\eta_{\eps}^+)''(u-v)
\sum_{k=1}^K \sum_{i,j=1}^d
\pxi \zeta_{ik}^a(u) \pyj \zeta_{jk}^a(v)
\varphi \dx \dt\dy\ds.
$$
Arguing exactly as in \cite{Bendahmane:2004go}, it follows that
\begin{equation}\label{eq:diff-diss-final}
	\begin{split}
		& \lim_{\eps\to 0} \Bigl( I_{\mathrm{diff}}(\eps)
		-\widetilde I_{\mathrm{diss}}(\eps)\Bigr)
		\\ & \quad
		\le \fourint \sum_{i,j=1}^d  r_{ij}^+(u,v)
		(\pxixj+2\pxiyj+\pyiyj) \varphi \dx\dt\dy\ds.
	\end{split}
\end{equation}

Fix a small number $\kappa>0$, and let us split $\Levy$ into two parts
\begin{align*}
	\Levy[\phi] & =\int_{\abs{z}\le \kappa}\left[ \phi(t,x+ z) - \phi(t,x)
	- z \cdot \nabla \phi \mathbf{1}_{|z|<1} \right]\, \pi(dz)
	\\ & \qquad +\int_{\abs{z}> \kappa}
	\left[ \phi(t,x+ z) - \phi(t,x) - z \cdot \nabla
	\phi \mathbf{1}_{|z|<1} \right]\, \pi(dz)
	\\ & =: \Levy_\kappa[\phi]+\Levy^\kappa[\phi],
	\qquad \forall \phi\in C^2,
\end{align*}
and similarly
$$
\Levy_x = \Levy_{x,\kappa} + \Levy_x^\kappa, \quad
\Levy_y = \Levy_{y,\kappa} + \Levy_y^\kappa, \quad
\Levy_{x+y} = \Levy_{x+y,\kappa} + \Levy_{x+y}^\kappa.
$$
The corresponding splitting of
$I_{\mathrm{fdiff}}(\eps)$ is written
$$
I_{\mathrm{fdiff}}(\eps)=
I_{\mathrm{fdiff},\kappa}(\eps)
+I_{\mathrm{fdiff}}^\kappa(\eps).
$$
We also need to introduce the
operator $\widetilde \Levy^\kappa$ defined by writing
\begin{align*}
	\Levy^\kappa[\varphi]
	& = \widetilde \Levy^\kappa [\varphi]
	-\left( \int_{\abs{z}>\kappa} z
	\mathbf{1}_{|z|<1}\,
	\pi(dz)\right) \cdot \nabla_x \varphi,
\end{align*}
with similar definitions for $\widetilde \Levy^\kappa_x$, $\widetilde \Levy^\kappa_y$,  and
$\widetilde \Levy^\kappa_{x+y}$. Observe that \eqref{eq:IBP}
continues to hold for all these operators. The function obtained
by replacing $\Levy^\kappa$ with $\widetilde \Levy^\kappa$ in the
definition of $I_{\mathrm{fdiff}}^\kappa(\eps)$ will be named
$\widetilde I_{\mathrm{fdiff}}^\kappa(\eps)$.

Clearly, in view of \eqref{intcond1},
\begin{equation}\label{eq:kappa-est}
	\abs{I_{\mathrm{fdiff},\kappa}(\eps)}\le
	C \norm{D^2\varphi}_{L^1(Q_T\times Q_T)}
	\int_{\abs{z}\le \kappa} \abs{z}^2\, \pi(dz)
	\overset{\kappa\to 0}{\longrightarrow} 0,
\end{equation}
for some constant $C$ independent of $\kappa$ and $\eps$.

Let us analyze $\widetilde I_{\mathrm{fdiff}}^\kappa(\eps)$.
By \eqref{eq:IBP},
\begin{align*}
	\widetilde I_{\mathrm{fdiff}}^\kappa(\eps)
	= \fourint  \Bigl(\widetilde \Levy^\kappa_x\left[\eta_\eps^+(u-v)\right]
	+ \widetilde \Levy^\kappa_y\left[\eta_\eps^+(u-v)\right]\Bigr)
	\varphi \dt \dx \dy\ds.
\end{align*}

Specifying $a=u(t,x)-v(s,y)$ and $b=u(t,x+z)-v(s,y)$ in \eqref{eq:Taylor} yields
\begin{equation}\label{eq:Taylor-1}
	\begin{split}
		&\eta_{\eps}^+(u(t,x+z)-v(s,y))-\eta_{\eps}^+(u(t,x)-v(s,y))
		\\ & = (\eta_{\eps}^+)'(u(t,x)-v(s,y)) \left(u(t,x+z) - u(t,x)\right)
		\\ & \qquad \qquad
		+\overline{(\eta_{\eps}^+)''}(u(t,\cdot)- v;z)\left(u(t,x+z) - u(t,x)\right)^2.
	\end{split}
\end{equation}

Similarly, taking $a=u(t,x)-v(s,y)$, $b=u(t,x)-v(s,y+z)$ in \eqref{eq:Taylor} yields
\begin{equation}\label{eq:Taylor-2}
	\begin{split}
		&\eta_{\eps}^+(u(t,x)-v(s,y+z))-\eta_{\eps}^+(u(t,x)-v(s,y))
		\\ & = -(\eta_{\eps}^+)'(u(t,x)-v(s,y)) \left(v(s,y+z) - v(s,y)\right)
		\\ & \qquad \qquad
		+\overline{(\eta_{\eps}^+)''}(u-v(s,\cdot);z) \left(v(s,y+z) - v(s,y)\right)^2.
	\end{split}
\end{equation}

Adding the first term on the right-hand side of \eqref{eq:Taylor-1} to the first
term on the right-hand side of \eqref{eq:Taylor-2} yields
\begin{align*}
	&(\eta_{\eps}^+)'(u(t,x)-v(s,y)) \left(u(t,x+z) - u(t,x)\right)\\
	& \qquad -(\eta_{\eps}^+)'(u(t,x)-v(s,y)) \left(v(s,y+z) - v(s,y)\right)
	\\ & = (\eta_{\eps}^+)'(u(t,x)-v(s,y))\Bigl[ \left(u(t,x+z) - v(s,y+z)\right)
	-\left(u(t,x) - v(s,y)\right)\Bigr]
	\\ & \le  \eta_\eps^+(u(t,x+z) - v(s,y+z))
	-\eta_\eps^+(u(t,x) - v(s,y)),
\end{align*}
where we have used the convexity of $\eta_\eps$ to derive the last inequality.

In view of these findings, we can rewrite 
$\widetilde I_{\mathrm{fdiff}}^\kappa(\eps)$ as follows:
\begin{align*}
	\widetilde I_{\mathrm{fdiff}}^\kappa(\eps)-I_{\mathrm{fdiss}}^\kappa(\eps)
	& \le  \fourint  \widetilde \Levy^\kappa_{x+y}
	\left[\eta_\eps^+(u(t,\cdot)-v(s,\cdot))\right]\varphi \dt \dx \dy\ds
	\\ & \overset{\text{\eqref{eq:IBP}}}{=}
	\fourint  \eta_\eps^+(u-v)\widetilde \Levy^\kappa_{x+y}[\varphi] \dt \dx \dy\ds,
\end{align*}
where
\begin{align*}
	I_{\mathrm{fdiss}}^\kappa(\eps) & =  \fourint\!\!\int_{\abs{z}>\kappa}
	\Biggl [ \overline{(\eta_{\eps}^+)''}(u(t,\cdot)-v;z) \left(u(t,x+z)-u(t,x)\right)^2
	\\ & \qquad\qquad\qquad\qquad\qquad
	+ \overline{(\eta_{\eps}^+)''}(u-v(s,\cdot);z) \left(v(s,y+z)-v(s,y)\right)^2 \Biggr]
	\\ & \qquad\qquad\qquad\qquad\qquad\qquad\qquad\qquad\qquad
	\times \varphi \, \pi(dz) \dx\dt\dy\ds.
\end{align*}
Consequently,
\begin{align*}
	I_{\mathrm{fdiff}}^\kappa(\eps)-I_{\mathrm{fdiss}}^\kappa(\eps)
	& \le \fourint  \eta_\eps^+(u-v)
	\Levy^\kappa_{x+y}[\varphi] \dt \dx \dy\ds,
\end{align*}

The next step is to first send $\kappa\to 0$ and then $\eps\to 0$. Related to this, observe that
$$
\lim_{\kappa\to 0} I_{\mathrm{fdiff}}^\kappa(\eps)= I_{\mathrm{fdiff}}(\eps), \quad
\lim_{\kappa\to 0} I_{\mathrm{fdiss}}^\kappa(\eps)= I_{\mathrm{fdiss}}(\eps)
$$
for each fixed $\eps>0$, by the dominated convergence theorem. Moreover, we clearly have
$\lim_{\kappa\to 0} \Levy^\kappa_{x+y}[\varphi] = \Levy_{x+y}[\varphi]$.
In view of this and \eqref{eq:kappa-est}, we conclude that
\begin{equation}\label{eq:fdiff-fdiss-final}
	I_{\mathrm{fdiff}}(\eps)-I_{\mathrm{fdiss}}(\eps)
	\le \fourint  \eta_\eps^+(u-v) \Levy[\varphi] \dt \dx \dy\ds.
\end{equation}

By \eqref{eq:diff-diss-final} and \eqref{eq:fdiff-fdiss-final}, It follows
from \eqref{in6} and sending $\eps\to 0$ that
\begin{equation}\label{in7}
	\begin{split}
		& \fourint \Biggl( (u-v)^+ (\pt+\ps)\varphi
		+ \sum_{i=1}^d  q_{i}^+(u,v) (\pxi+\pyi) \varphi
		\\ & \qquad
		+\sum_{i,j=1}^d  r_{ij}^+(u,v)
		(\pxixj+2\pxiyj+\pyiyj) \varphi
		+ \eta^+(u-v) \Levy_{x+y}[\varphi] \Biggr)\dx\dt\dy\ds
		\\ & + \threeint (u_0-v)^+ \varphi(0,x,s,y) \dx\dy\ds
		+ \threeint (u-v_0)^+ \varphi(t,x,0,y) \dx\dt\dy \ge 0.
	\end{split}
\end{equation}

Let us specify the test function $\varphi=\varphi(t,x,s,y)$. To this end, fix a nonnegative
test function $\phi=\phi(t,x)\in C_c^{\infty}([0,\infty) \times \R^d)$, and
pick two sequences $\Set{\theta_{\nu}}_{\nu>0} \subset C_c^{\infty}(0,\nu)$,
$\Set{\delta_{\mu}}_{\mu>0} \subset C_c^{\infty}(B(0,\mu))$ of
approximate delta functions, where $B(0,\mu)$ denotes the open ball
centered at the origin with radius $\mu$. Then take
\begin{equation}\label{eq:testfunc}
	\varphi(t,x,s,y) = \theta_{\nu}(s-t) 
	\delta_{\mu}(y-x) \phi(t,x).
\end{equation}

Simple calculations reveal that
\begin{align*}
	(\pt + \ps) \varphi
	& =\theta_{\nu}(s-t) \delta_{\mu}(y-x) \pt \phi(t,x),
	\\ (\pxi + \pyi)\varphi
	& = \theta_{\nu}(s-t) \delta_{\mu}(y-x) \pxi \phi(t,x),
	\\ (\pxixj +2 \pxiyj + \pyiyj) \varphi
	& = \theta_{\nu}(s-t) \delta_{\mu}(y-x) \pxixj \phi(t,x)
\end{align*}
and
\begin{align*}
	&\varphi(t,x+z,s,y+z) - \varphi(t,x,s,y)
	\\ & \qquad
	= \theta_{\nu}(s-t) \delta_{\mu}(y-x)
	\left( \phi(t,x+z)-\phi(t,x) \right).
\end{align*}

Note that $\theta_{\nu} = 0 $ on $(-\infty,0]$ and so $\varphi(t,x,0,y)\equiv 0.$
By the choice of the test function $\varphi$ and the observations above, we
deduce from \eqref{in7} that
\begin{equation}\label{in8}
	\begin{split}
		& \fourint (u-v)^+ \theta_{\nu}(s-t) \delta_{\mu}(y-x)
		\pt \phi(t,x) \dx\dt\dy\ds \\
		&\quad +\fourint
		\sum_{i=1}^d q_i^+(u,v) \theta_{\nu}(s-t)
		\delta_{\mu}(y-x)\pxi   \phi(t,x)   \dxdt \dyds \\
		&\quad + \fourint \sum_{i,j=1}^d r_{ij}^+(u,v)
		\theta_{\nu}(s-t) \delta_{\mu}(y-x)\pxixj \phi(t,x) \dx\dt\dy\ds
		\\ & \quad + \fourint (u-v)^+
		\theta_{\nu}(s-t) \delta_{\mu}(y-x) \Levy[\phi] \dx\dt\dy\ds
		+ I_{u_0,v}(\nu,\mu)\geq 0,
	\end{split}
\end{equation}
where
\begin{align*}
	I_{u_0,v}(\nu,\mu) & :=
	\threeint (u_0-v)^+
	\theta_{\nu}(s)\delta_{\mu}(y-x) \phi(0,x) \dx\dy\ds
	\\ & = -\threeint (u_0-v)^+
	\ps\left(\tilde{\phi}_{\nu}(s)
	\delta_{\mu}(y-x)\phi(0,x)\right) \dx\dy\ds,
\end{align*}
with
$$
\tilde{\phi}_{\nu}(s) := \int_s^T \theta_{\nu}(\tau)\dtau
= \int_{\min(s,\nu)}^{\nu} \theta_{\nu}(\tau)\dtau
\overset{\nu\to 0}{\longrightarrow} 1.
$$
Specifying $\varphi=\tilde{\phi}_{\nu}(s)\delta_{\mu}(y-x)$ in
the entropy inequality for $v$ and noting
that $\theta_{\nu}(s)$ vanishes for $s>\nu$, we obtain
\begin{equation}\label{in9}
	\begin{split}
		& \iint  (u_0-v)^+\ps \varphi(s,x,y) \dy\ds \\ & \qquad
		\le \iint (u_0 -v)^+\theta_{\nu}(s) 
		\delta_{\mu}(y-x)\phi(0,x) \dy\ds+o(\nu)
		\\ & \quad \overset{\nu\to 0}{\longrightarrow}
		\iint (u_0 -v)^+ \delta_{\mu}(y-x) \phi(0,x)\dy\ds,
	\end{split}
\end{equation}
where the ``$o(\nu)$" term follows from an integrability argument.

Hence, sending $\nu,\mu \to 0,$ we deduce
\begin{equation}\label{in11}
	\begin{split}
		& \limsup_{\mu \to 0} \limsup_{\nu \to 0} I_{u_0,v}(\nu,\mu)
		\\ & \quad \leq \limsup_{\mu \to 0} \iint
		(u_0 -v_0)^+  \delta_{\mu}(y-x) \phi(0,x)\dx\dy
		\\ & \quad
		= \int (u_0-v_0)^+\phi(0,x)\dx,
	\end{split}
\end{equation}
with $u_0=u_0(x)$ and $v_0=v_0(x)$.

Keeping in mind \eqref{in11} when sending $\mu,\nu \to 0$ 
in \eqref{in8}, we conclude that
\begin{equation}\label{in12}
	\begin{split}
		& \iint_{Q_T} \Biggl( (u -v)^+ \pt \phi
		+ \sum_{i=1}^d q_i^+(u,v)  \pxi \phi
		\\ & \qquad\qquad\qquad\quad
		+ \sum_{i,j=1}^d r_{ij}^+(u,v) \pxixj \phi
		+(u-v)^+\Levy[\phi] \Biggr) \dx\dt
		\\ & \qquad\qquad\qquad\qquad\quad
		+ \int_{\R^d} (u_0-v_0)^+ \phi(0,x) \dx \geq 0,
    \end{split}
\end{equation}
where all the involved functions depend on $(t,x)$.
It now only takes a standard argument to conclude from
\eqref{in12} that Theorem \ref {th:unique} holds. Indeed, one chooses
a sequence of functions $0 \leq \phi \leq 1$ from $C^\infty_c([0,T)\times \R^d)$
that converges to $\mathbf{1}_{[0,t) \times \R^d}$ for a Lebesgue point $t$ of
$\int_{\R^d} (u -v)^+\dx$ and then use the integrability of $u,v$ to conclude the proof.

Finally, let us prove the last part of Theorem \ref{th:unique}, that is, we shall 
establish the $L^1$ contraction property for entropy solutions satisfying 
the simpler entropy condition \eqref{entineq-simpler} in which the 
fractional parabolic dissipation measure has been dropped. 
To this end, let us assume that \eqref{eq:add-cond} holds. 
Arguing exactly as before we arrive at 
\eqref{in6} without the $I_{\mathrm{fdiss}}(\eps)$ term:
\begin{equation}\label{intineq}
	\begin{split}
		& \fourint \Biggl( (u-v)^+ (\pt+\ps)\varphi
		+ \sum_{i=1}^d  q_{i}^+ (u,v) (\pxi+\pyi) \varphi
		\\ & \qquad\qquad\qquad\qquad
		+\sum_{i,j=1}^d  r_{ij}^+(u,v)
		(\pxixj+2\pxiyj+\pyiyj) \varphi
		\\ & \qquad \qquad \qquad\qquad\qquad
		+ (u-v)^+ (\Levy_{x}[\varphi]+\Levy_y[\varphi]) 
		\Biggr)\dx\dt\dy\ds
		\\ & \quad 
		+ \threeint (u_0-v)^+\varphi|_{t=0}\dx\dy\ds 
		+ \threeint (u-v_0)^+\varphi|_{s=0}\dx\dt\dy \ge 0. 
	\end{split}
\end{equation}

The new treatment concerns the fractional term 
in \eqref{intineq} only, which we denote by $J(\nu,\mu)$, 
i.e., $J(\nu,\mu):=\fourint(u-v)^+ (\Levy_{x}[\varphi]+\Levy_y[\varphi]) \dx\dt\dy\ds$; we 
employ the same test function $\varphi$ as before, cf.~\eqref{eq:testfunc}. 
By letting $z \mapsto -z$ 
in the $\Levy_y$ term, keeping in mind 
that $\pi(d(-z))=-\pi(dz)$, we obtain
\begin{equation}\label{fracinter}
	\begin{split}
		J(\nu,\mu)=\fourint \int_{\R^d \setminus \{0\}} 
		& (u(t,x)-v(s,y))^+ 
		\theta_{\nu}(s-t) 
		\\ & \times
		\Bigg(\delta_{\mu}(y-x-z) (\phi(t,x+z)-\phi(t,x)) 
		\\ & \qquad\qquad\quad
		-\delta_{\mu}(y-x) \nabla \phi(t,x) 
		\cdot z\, 1_{|z|<1} \Bigg)\, \pi(dz) \dx\dt\dy\ds
	\end{split}
\end{equation}
Sending the ``temporal parameter" $\nu$ to zero in \eqref{fracinter} yields
\begin{align*}
		J(\mu) & :=\lim_{\nu\to 0} J(\nu,\mu)
		\\ & =\fourint \int_{\R^d \setminus \{0\}}
		(u(t,x)-v(s,y))^+ 
		\theta_{\nu}(s-t) 
		\\ & \qquad\qquad\qquad\quad \times
		\Bigg(\delta_{\mu}(y-x-z) (\phi(t,x+z)-\phi(t,x)) 
		\\ & \qquad\qquad\quad\qquad\qquad\qquad\quad
		-\delta_{\mu}(y-x) \nabla \phi(t,x) 
		\cdot z\, 1_{|z|<1} \Bigg)\,\pi(dz) \dx\dt\dy\ds
\end{align*}

Under the additional requirements listed in \eqref{eq:add-cond} 
we are also allowed to send the ``spatial parameter" $\mu$ to zero, resulting in
\begin{equation}\label{Jfin1}
	\begin{split}
		J & := \lim_{\mu\to 0} J(\mu) 
		\\ & = \iint \int_{\R^d \setminus \{0\}} 
		\Biggl( (u(t,x) -v(t,x+z))^+ 
		\left( \phi(t,x+z) - \phi(t,x) \right) \\
		&\qquad\qquad\qquad\qquad\quad
		-(u(t,x)-v(t,x))^+ \nabla \phi(t,x) 
		\cdot z\, 1_{|z|<1} \Biggr) \,\pi(dz) \dx\dt.
	\end{split}
\end{equation}
We divide the remaining discussion into two cases.\\

\textbf{Case 1: $|z|\pi(dz) \in L^1(\R^d \setminus \{0\})$.} 
Adding and subtracting identical terms we obtain
\begin{equation}\label{J}
		J = \iint (u(t,x)-v(t,x))^+
		\Levy[\phi](t,x) \dt \dx + E,
\end{equation}
where		 
\begin{align*}
	|E| &\le \iint \int_{\R^d \setminus \{0\}} 
	|v(t,x+z)-v(t,x)|\, |\phi(t,x+z)-\phi(t,x)| \,\pi(dz) \dt \dx
	\\ & \leq \Bigl(2 T\norm{v}_{L^\infty(0,T;L^1(\R^d))} 
	\norm{\phi'}_{L^\infty(Q_T)}\Bigr) |z| 
	\in L^1(\R^d \setminus \{0\};\pi(dz)),
\end{align*}
so that we can employ the dominated convergence theorem 
to send $\phi \to 1$, and consequently $|E|\to 0$.\\

\textbf{Case 2: $|z|^2 \pi(dz) \in L^1(\R^d \setminus \{0\})$ 
and $v \in L^\infty(0,T;BV(\R^d))$.} In this case the error 
term $|E|$  in (\ref{J}) can be estimated as follows:
$$  
|E|\leq \Bigl( T \norm{v}_{L^\infty(0,T;BV(\R^d)}
\norm{\phi''}_{L^\infty(Q_T)}\Bigr) 
|z|^2 \in L^1(\R^d \setminus \{0\};\pi(dz)),
$$
and again we can employ the dominated convergence theorem. 

This concludes the proof of Theorem \ref{th:unique}.

\section{Proof of Theorem \ref{th:contdep} (continuous dependence)}\label{prcdep}
We again employ the doubling of variables 
device as in the previous section, but with a slightly different 
choice of the entropy function. For each $\eps>0$, define
$$
\sgn_{\eps}(\xi) =
\begin{cases}
	-1, & \textrm{if $\xi<-\eps$}\\
	\sin(\frac{\pi}{2\eps}\xi), & \textrm{if $ |\xi| \leq \eps$}\\
	1, & \textrm{if $\xi > \eps$},
\end{cases}
$$
which is a $C^1$  approximation of $\sgn(\cdot)$.
This choice gives rise to a $C^2$ approximation
$\eta_\eps(z)=\int_0^z \sgn_{\eps}(\xi)\,d\xi$ of the entropy flux $\abs{z}$.
As before, we introduce the corresponding entropy flux functions
$\eta^{\eps}(u,c)$, $q_i^{\eps}(u,c)$, and $r_{ij}^{\eps}(u,c)$.
We now employ the doubling variables technique using the test function
$$
\varphi(t,x,s,y) =   \theta_{\nu}(s-t) \delta_{\mu}(y-x)  \Theta_{\alpha}(t),
$$
where $\theta_{\nu}$, $\delta_{\mu}$ are symmetric approximate
delta functions with support in $(-\nu,\nu)$ and $B(0,\mu)$, respectively.
Fix a time $\tau$ from $(0,T)$.
For any $\alpha>0$ with $0<\alpha<\min(\tau_0,T-\tau)$, we define
$$
\Theta_{\alpha}(t)
=H_{\alpha}(t) - H_{\alpha}(t-\tau),
\quad
H_{\alpha}(t)=\int_{-\infty}^{t}
\theta_{\alpha}(\sigma)\,d\sigma.
$$
so that $\Theta_{\alpha}'(t)=\theta_{\alpha}(t)-\theta_{\alpha}(t-\tau)$.

Proceeding as in the previous section (cf.~also \cite{Chen:2005wf}) and
sending $\eps \to 0$, we find
$$
-\fourint
\abs{u-v} \theta_{\nu}(s-t) \delta_{\mu}(y-x)
\Theta_\alpha'(t)\dx\dt\dy\ds
\le I_{\mathrm{conv}} - I_{\mathrm{diff}}+I_{\mathrm{fdiff}},
$$
where
$$
I_{\mathrm{conv}} := \fourint \left[ G(u,v) - F(u,v) \right]\cdot
\nabla_x \delta_{\mu}(y-x) \theta_{\nu}(s-t) \Theta_\alpha(t)\dx\dt\dy\ds,
$$
$$
F(u,v) := \sgn(u-v) \left(f(u) - f(v)\right), \quad
G(u,v) := \sgn(u-v) \left(g(u) - g(v)\right),
$$
$$
I_{\mathrm{diff}} :=
\fourint \sum_{i,j=1}^d \Theta_\alpha(t) \theta_{\nu}(s-t)
\pxixj\delta_{\mu}(y-x) \int_v^u \sgn(\xi -v)
\eps_{ij}^{a-b}(\xi) \dxi \dx\dt\dy\ds,
$$
$$
\eps_{ij}^{a-b}(\xi) := \sum_{k=1}^K \left(\sigma_{ik}^a(\xi) \sigma_{jk}^a(\xi)
- 2 \sigma_{ik}^a(\xi) \sigma_{jk}^b(\xi) + \sigma_{ik}^b(\xi)\sigma_{jk}^b(\xi)\right).
$$
and $I_{\mathrm{fdiff}}=
I_{\mathrm{fdiff}_1}+I_{\mathrm{fdiff}_2}$ with 
\begin{align*}
	I_{\mathrm{fdiff}_{1}} :=
	\fourint \int_{|z|<1} &\abs{u-v} 
	\theta_{\nu}(s-t)\Theta_{\alpha}(t)
	\\ & \quad \times 
	\Bigl[\delta_{\mu}(y-x-z)-\delta_{\mu}(y-x)
	-\nabla \delta_{\mu}(y-x) \cdot z \Bigr] 
	\\ & \quad \qquad \times (m(z) - \tilde{m}(z)) 
	\dz\dx\dt\dy\ds
\end{align*}
and
\begin{align*}
	I_{\mathrm{fdiff}_{2}}:=
	\fourint \int_{|z|\geq 1} 
	\abs{u-v}&\theta_{\nu}(s-t)\Theta_{\alpha}(t)
	\Bigl[\delta_{\mu}(y-x-z)-\delta_{\mu}(y-x) \Bigr]
	\\ & \times (m(z)-\tilde{m}(z)) 
	\dz\dx\dt\dy\ds,
\end{align*}

By triangle inequality
\begin{align*}
	& -\fourint \abs{u(t,x)-v(s,y)} \theta_{\nu}(s-t) \delta_{\mu}(y-x)
	\Theta_\alpha'(t)\dx\dt\dy\ds
	\\ & \qquad \geq -\fourint \abs{u(t,y)-v(t,y)} \theta_{\nu}(s-t) \delta_{\mu}(y-x)
	\abs{\Theta_\alpha'(t)}\dx\dt\dy\ds
	\\ & \quad \qquad - \fourint \abs{v(t,y)-v(s,y)} \theta_{\nu}(s-t) \delta_{\mu}(y-x)
	\abs{\Theta_\alpha'(t)}\dx\dt\dy\ds
	\\ & \quad \qquad - \fourint \abs{u(t,x)-u(t,y)} \theta_{\nu}(s-t) \delta_{\mu}(y-x)
	\abs{\Theta_\alpha'(t)}\dx\dt\dy\ds
	\\ & =: L + R_t + R_x.
\end{align*}

Keeping in mind that $v\in C(L^1)$ and $u\in L^\infty(BV)$, it
is standard to show that
$$
\lim_{\nu \to 0} R_t = 0, \quad
\limsup_{\alpha\to 0}\abs{R_x} \leq C\mu
$$
and moreover, since also $u(t)\to u_0,v(t)\to v_0$ as $t\to 0$,
$$
\lim_{\alpha\to 0} L
=  \norm{u(\tau,\cdot) - v(\tau,\cdot)}_{L^1(\R^d)}
-\norm{u_0-v_0}_{L^1(\R^d)}.
$$

Following \cite{Chen:2005wf}, using $u\in L^\infty(BV)$ we conclude that
$$
\lim_{\alpha \to 0} \lim_{\nu \to 0} \abs{I_{\mathrm{conv}}}
\leq C \tau \norm{f-g}_{\mathrm{Lip}(I)},
$$
and, exploiting also that $\int \abs{\pxi \delta_\mu}\le C/\mu$,
$$
\lim_{\alpha \to 0} \lim_{\nu \to 0} \abs{I_{\mathrm{diff}}}
\leq \frac{C}{\mu}
\tau \norm{(\sigma^a - \sigma^b)
(\sigma^a -\sigma^b)^{\mathrm{tr}}}_{L^{\infty}(I;\R^{d\times d})}.
$$

It remains to estimate $\abs{I_{\mathrm{fdiff}}}$. 
First, we nconsider $I_{\mathrm{fdiff}_{1}}$. 
Using the Taylor and Fubini theorems we obtain
\begin{align*}
	\abs{I_{\mathrm{fdiff}_{1}}} 
	&=\threeint \int_{|z|<1} \int_0^1  (1-\tau) 
	\theta_{\mu}(s-t) \Theta_{\alpha}(t) (\tilde{m}(z)-m(z))
	\\ & \qquad \times 
	\left( \int_{\R^d} \abs{u(t,x)-v(s,y)} D^2 
	\delta_{\mu}(y-x-\tau z)\, z \cdot z
	\dx \right)\dtau \dz \dy \ds \dt.
\end{align*}
Thanks to $|u(t,\cdot)-v(s,y)| \in BV(\R^d)$, an integration by parts yields
\begin{equation}\label{i22bdd}
	\begin{split}
		I_{\mathrm{fdiff}_{1}} &=
		\threeint \int_{|z|<1} \int_0^1  (1-\tau) 
		\theta_{\mu}(s-t) \Theta_{\alpha}(t) (\tilde{m}(z)-m(z)) \\
		&\times \left( \int_{\R^d} \nabla \delta_{\mu}(y-x-\tau z)\cdot z\,
		D_x\left(\abs{u(t,x)-v(s,y)}\right)\cdot z \dx \right) 
		\dtau \dz\dy\ds\dt,
	\end{split}
\end{equation}
where the inner integral is taken with respect to the bounded Borel 
measure $D\left(|u(t,\cdot)-v(s,y)|\right)\cdot z$. 
Since $|D(u(t,\cdot) - v(s,y))| \leq |D(u(t,\cdot))|$, the term inside the 
parentheses in \eqref{i22bdd}, is upper bounded by
\begin{align*}
	|z|^2 \int_{\R^d} \int_{\R^d} 
	|\nabla \delta_{\mu}(y-x-\tau z)|\,|dD(u(t,\cdot))(x)|\dy
	\le |z|^2\abs{u(t,\cdot)}_{BV(\R^d)} \norm{\nabla \delta_{\mu}}_{L^1(\R^d)},
\end{align*}
where we have used that $|D u(t,\cdot)|$ is finite and the Fubini's theorem to 
first integrate with respect to $y$. Hence, 
$$
\lim_{\alpha \to 0} \lim_{\nu \to 0} 
\abs{I_{\mathrm{fdiff}_{1}}} 
\leq \frac{C}{\mu} 
\tau \int_{|z|<1} |z|^2 
\abs{m(z)-\tilde{m}(z)}\dz,
$$
where $C>0$ is a finite constant.

Similarly, relying again on the $L^\infty(BV)$ regularity of $u$, it 
is not difficult to deduce via an integration by parts the estimate
$$
\lim_{\nu \to 0} \lim_{\alpha \to 0} 
\abs{I_{\mathrm{fdiff}_{2}}}
\leq C\tau\int_{|z|\geq 1}
|z|\abs{m(z)-\tilde{m}(z)}\dz.
$$

Finally, we collect the bounds we have obtained so far and then 
optimize over $\mu$ to obtain the desired 
continuous dependence estimate \eqref{contdest}.


\begin{thebibliography}{10}

\bibitem{Alibaud:2007mi}
N.~Alibaud.
\newblock Entropy formulation for fractal conservation laws.
\newblock {\em J. Evol. Equ.}, 7(1):145--175, 2007.

\bibitem{Alibaud:2007dw}
N.~Alibaud.
\newblock Existence, uniqueness and regularity for nonlinear parabolic
  equations with nonlocal terms equations with nonlocal terms.
\newblock {\em NoDEA Nonlinear Differential Equations Appl.}, 14(3-4):259--289,
  2007.

\bibitem{Alibaud:2007qe}
N.~Alibaud, J.~Droniou, and J.~Vovelle.
\newblock Occurrence and non-appearance of shocks in fractal {B}urgers
  equations.
\newblock {\em J. Hyperbolic Differ. Equ.}, 4(3):479--499, 2007.

\bibitem{Alvarez:1996lq}
O.~Alvarez and A.~Tourin.
\newblock Viscosity solutions of nonlinear integro-differential equations.
\newblock {\em Ann. Inst. H. Poincar\'e Anal. Non Lin\'eaire}, 13(3):293--317,
  1996.

\bibitem{Arisawa:2006db}
M.~Arisawa.
\newblock A new definition of viscosity solutions for a class of second-order
  degenerate elliptic integro-differential equations.
\newblock {\em Ann. Inst. H. Poincar\'e Anal. Non Lin\'eaire}, 23(5):695--711,
  2006.

\bibitem{Barles:1997xj}
G.~Barles, R.~Buckdahn, and E.~Pardoux.
\newblock Backward stochastic differential equations and integral-partial
  differential equations.
\newblock {\em Stochastics Stochastics Rep.}, 60(1-2):57--83, 1997.

\bibitem{Barles:2008lq}
G.~Barles and C.~Imbert.
\newblock Second-order elliptic integro-differential equations: viscosity
  solutions' theory revisited.
\newblock {\em Ann. Inst. H. Poincar\'e Anal. Non Lin\'eaire}, 25(3):567--585,
  2008.

\bibitem{Bendahmane:2004go}
M.~Bendahmane and K.~H. Karlsen.
\newblock Renormalized entropy solutions for quasi-linear anisotropic
  degenerate parabolic equations.
\newblock {\em SIAM J. Math. Anal.}, 36(2):405--422 (electronic), 2004.

\bibitem{Benth:2001tv}
F.~E. Benth, K.~H. Karlsen, and K.~Reikvam.
\newblock Optimal portfolio selection with consumption and nonlinear
  integrodifferential equations with gradient constraint: a viscosity solution
  approach.
\newblock {\em Finance and Stochastic}, 5:275--303, 2001.

\bibitem{Benth:2002mk}
F.~E. Benth, K.~H. Karlsen, and K.~Reikvam.
\newblock Portfolio optimization in a {L}\'evy market with intertemporal
  substitution and transaction costs.
\newblock {\em Stoch. Stoch. Rep.}, 74(3-4):517--569, 2002.

\bibitem{Bertoin:1996ud}
J.~Bertoin.
\newblock {\em L\'evy processes}, volume 121 of {\em Cambridge Tracts in
  Mathematics}.
\newblock Cambridge University Press, Cambridge, 1996.

\bibitem{Biler:1998ai}
P.~Biler, T.~Funaki, and W.~A. Woyczynski.
\newblock Fractal {B}urgers equations.
\newblock {\em J. Differential Equations}, 148(1):9--46, 1998.

\bibitem{Biler:2001th}
P.~Biler, G.~Karch, and W.~A. Woyczy{\'n}ski.
\newblock Asymptotics for conservation laws involving {L}\'evy diffusion
  generators.
\newblock {\em Studia Math.}, 148(2):171--192, 2001.

\bibitem{Bossy:2002eu}
M.~Bossy and B.~Jourdain.
\newblock Rate of convergence of a particle method for the solution of a 1{D}
  viscous scalar conservation law in a bounded interval.
\newblock {\em Ann. Probab.}, 30(4):1797--1832, 2002.

\bibitem{Brandolese:2008sp}
L.~Brandolese and G.~Karch.
\newblock Far field asymptotics of solutions to convection equation with
  anomalous diffusion.
\newblock {\em J. Evol. Equ.}, 8(2):307--326, 2008.

\bibitem{Caffarelli:2007hh}
L.~Caffarelli and L.~Silvestre.
\newblock An extension problem related to the fractional {L}aplacian.
\newblock {\em Comm. Partial Differential Equations}, 32(7-9):1245--1260, 2007.

\bibitem{Caffarelli:2007qr}
L.~Caffarelli and L.~Silvestre.
\newblock Regularity theory for fully nonlinear integro-differential equations.
\newblock Submitted, 2007.

\bibitem{Caffarelli:2008bx}
L.~A. Caffarelli, S.~Salsa, and L.~Silvestre.
\newblock Regularity estimates for the solution and the free boundary of the
  obstacle problem for the fractional {L}aplacian.
\newblock {\em Invent. Math.}, 171(2):425--461, 2008.

\bibitem{Carrillo:1999hq}
J.~Carrillo.
\newblock Entropy solutions for nonlinear degenerate problems.
\newblock {\em Arch. Ration. Mech. Anal.}, 147(4):269--361, 1999.

\bibitem{Chen:2005wf}
G.-Q. Chen and K.~H. Karlsen.
\newblock Quasilinear anisotropic degenerate parabolic equations with
  time-space dependent diffusion coefficients.
\newblock {\em Commun. Pure Appl. Anal.}, 4(2):241--266, 2005.

\bibitem{Chen:2006oy}
G.-Q. Chen and K.~H. Karlsen.
\newblock {$L\sp 1$}-framework for continuous dependence and error estimates
  for quasilinear anisotropic degenerate parabolic equations.
\newblock {\em Trans. Amer. Math. Soc.}, 358(3):937--963 (electronic), 2006.

\bibitem{Chen:2003td}
G.-Q. Chen and B.~Perthame.
\newblock Well-posedness for non-isotropic degenerate parabolic-hyperbolic
  equations.
\newblock {\em Ann. Inst. H. Poincar\'e Anal. Non Lin\'eaire}, 20(4):645--668,
  2003.

\bibitem{Cont:2004gk}
R.~Cont and P.~Tankov.
\newblock {\em Financial modelling with jump processes}.
\newblock Chapman \& Hall/CRC Financial Mathematics Series. Chapman \&
  Hall/CRC, Boca Raton, FL, 2004.

\bibitem{Droniou:2003mz}
J.~Droniou, T.~Gallouet, and J.~Vovelle.
\newblock Global solution and smoothing effect for a non-local regularization
  of a hyperbolic equation.
\newblock {\em J. Evol. Equ.}, 3(3):499--521, 2003.
\newblock Dedicated to Philippe B{\'e}nilan.

\bibitem{Droniou:2006os}
J.~Droniou and C.~Imbert.
\newblock Fractal first-order partial differential equations.
\newblock {\em Arch. Ration. Mech. Anal.}, 182(2):299--331, 2006.

\bibitem{Ev}
L. C.~Evans.
\newblock Partial Differential Equations.
\newblock {\em Providence, RI: American Mathematical Society}, 2002.

\bibitem{Garroni:2002il}
M.~G. Garroni and J.~L. Menaldi.
\newblock {\em Second order elliptic integro-differential problems}, volume 430
  of {\em Chapman \& Hall/CRC Research Notes in Mathematics}.
\newblock Chapman \& Hall/CRC, Boca Raton, FL, 2002.

\bibitem{Jacob:2001rf}
N.~Jacob.
\newblock {\em Pseudo differential operators and {M}arkov processes. {V}ol.
  {I}}.
\newblock Imperial College Press, London, 2001.
\newblock Fourier analysis and semigroups.

\bibitem{Jacob:2002xp}
N.~Jacob.
\newblock {\em Pseudo differential operators \& {M}arkov processes. {V}ol.
  {II}}.
\newblock Imperial College Press, London, 2002.
\newblock Generators and their potential theory.

\bibitem{Jacob:2005ss}
N.~Jacob.
\newblock {\em Pseudo differential operators and {M}arkov processes. {V}ol.
  {III}}.
\newblock Imperial College Press, London, 2005.
\newblock Markov processes and applications.

\bibitem{Jakobsen:2005jy}
E.~R. Jakobsen and K.~H. Karlsen.
\newblock Continuous dependence estimates for viscosity solutions of
  integro-{PDE}s.
\newblock {\em J. Differential Equations}, 212(2):278--318, 2005.

\bibitem{Jakobsen:2006aa}
E.~R. Jakobsen and K.~H. Karlsen.
\newblock A "maximum principle for semicontinuous functions" applicable to
  integro-partial differential equations.
\newblock {\em NoDEA Nonlinear Differential Equations Appl.}, 13(2):137--165,
  2006.

\bibitem{Karch:2008dp}
G.~Karch, C.~Miao, and X.~Xu.
\newblock On convergence of solutions of fractal {B}urgers equation toward
  rarefaction waves.
\newblock {\em SIAM J. Math. Anal.}, 39(5):1536--1549, 2008.

\bibitem{Karlsen:2002bh}
K.~H. Karlsen and M.~Ohlberger.
\newblock A note on the uniqueness of entropy solutions of nonlinear degenerate
  parabolic equations.
\newblock {\em J. Math. Anal. Appl.}, 275(1):439--458, 2002.

\bibitem{Karlsen:2003za}
K.~H. Karlsen and N.~H. Risebro.
\newblock On the uniqueness and stability of entropy solutions of nonlinear
  degenerate parabolic equations with rough coefficients.
\newblock {\em Discrete Contin. Dyn. Syst.}, 9(5):1081--1104, 2003.

\bibitem{KU2}
K.H. ~Karlsen and S. ~Ulusoy.
\newblock Existence and numerics for entropy solutions to 
fractional quasilinear anisotropic degenerate parabolic equations. in progress.

\bibitem{Kiselev:2008jt}
A.~Kiselev, F.~Nazarov, and R.~Shterenberg.
\newblock Blow up and regularity for fractal {B}urgers equation.
\newblock Submitted, 2008.

\bibitem{Kruzkov:1970kx}
S.~N. Kru{\v{z}}kov.
\newblock First order quasilinear equations with several independent variables.
\newblock {\em Mat. Sb. (N.S.)}, 81 (123):228--255, 1970.

\bibitem{Mascia:2002dq}
C.~Mascia, A.~Porretta, and A.~Terracina.
\newblock Nonhomogeneous {D}irichlet problems for degenerate
  parabolic-hyperbolic equations.
\newblock {\em Arch. Ration. Mech. Anal.}, 163(2):87--124, 2002.

\bibitem{Michel:2003tw}
A.~Michel and J.~Vovelle.
\newblock Entropy formulation for parabolic degenerate equations with general
  {D}irichlet boundary conditions and application to the convergence of {FV}
  methods.
\newblock {\em SIAM J. Numer. Anal.}, 41(6):2262--2293 (electronic), 2003.

\bibitem{Perthame:2003yq}
B.~Perthame and P.~E. Souganidis.
\newblock Dissipative and entropy solutions to non-isotropic degenerate
  parabolic balance laws.
\newblock {\em Arch. Ration. Mech. Anal.}, 170(4):359--370, 2003.

\bibitem{Pham:1998zt}
H.~Pham.
\newblock Optimal stopping of controlled jump diffusion processes: a viscosity
  solution approach.
\newblock {\em J. Math. Systems Estim. Control}, 8(1):27 pp.\ (electronic),
  1998.

\bibitem{Sato:1999qd}
K.-i. Sato.
\newblock {\em L\'evy processes and infinitely divisible distributions},
  volume~68 of {\em Cambridge Studies in Advanced Mathematics}.
\newblock Cambridge University Press, Cambridge, 1999.

\bibitem{Sayah:1991xy}
A.~Sayah.
\newblock \'{E}quations d'{H}amilton-{J}acobi du premier ordre avec termes
  int\'egro-diff\'erentiels. {I}. {U}nicit\'e des solutions de viscosit\'e.
\newblock {\em Comm. Partial Differential Equations}, 16(6-7):1057--1074, 1991.

\bibitem{Sayah:1991jk}
A.~Sayah.
\newblock \'{E}quations d'{H}amilton-{J}acobi du premier ordre avec termes
  int\'egro-diff\'erentiels. {II}. {E}xistence de solutions de viscosit\'e.
\newblock {\em Comm. Partial Differential Equations}, 16(6-7):1075--1093, 1991.

\bibitem{Silvestre:2006bq}
L.~Silvestre.
\newblock H\"older estimates for solutions of integro-differential equations
  like the fractional {L}aplace.
\newblock {\em Indiana Univ. Math. J.}, 55(3):1155--1174, 2006.

\bibitem{Silvestre:2007qp}
L.~Silvestre.
\newblock Regularity of the obstacle problem for a fractional power of the
  {L}aplace operator.
\newblock {\em Comm. Pure Appl. Math.}, 60(1):67--112, 2007.

\bibitem{Soner:IntegroDif}
H.~M. Soner.
\newblock Optimal control of jump-{M}arkov processes and viscosity solutions.
\newblock In {\em Stochastic differential systems, stochastic control theory
  and applications (Minneapolis, Minn., 1986)}, volume~10 of {\em IMA Vol.
  Math. Appl.}, pages 501--511. Springer, New York, 1988.

\end{thebibliography}
\end{document}